\documentclass[11pt]{article}
\usepackage{amssymb}
\usepackage{amsthm}
\usepackage{amsmath,amsfonts,amssymb}

\usepackage{graphicx}
\usepackage{algpseudocode}
\usepackage{algorithm}
\usepackage{natbib}
\usepackage{color}
\usepackage{multirow}

\def\bel{\begin{equation}\label}
\def\eeq{\end{equation}}

\begin{document}

\title{Day-to-day dynamic traffic assignment model with variable message signs and endogenous user compliance}
\author{Ke Han$\thanks{Corresponding author, e-mail: k.han@imperial.ac.uk;}$
\qquad Margherita Mascia$\thanks{e-mail: m.mascia@imperial.ac.uk;}$
\qquad Robin North$\thanks{e-mail: robin.north@imperial.ac.uk;}$
\qquad Simon Hu$\thanks{e-mail: j.s.hu05@imperial.ac.uk;}$ 
\qquad Gabriel Eve$\thanks{e-mail: gabriel.eve09@imperial.ac.uk;}$\\\\
\textit{Center for Transport Studies}\\
\textit{Department of Civil and Environmental Engineering}\\
\textit{Imperial College London, SW7 2BU, UK}\\}
\date{\today}
\maketitle

\begin{abstract}
This paper proposes a dual-time-scale, day-to-day dynamic traffic assignment model that takes into account {\it variable message signs} (VMS) and its interactions with drivers' travel choices and adaptive learning processes. The within-day dynamic is captured by a {\it dynamic network loading} problem with en route update of path choices influenced by the VMS; the day-to-day dynamic is captured by a simultaneous route-and-departure-time adjustment process that employs bounded user rationality. Moreover, we describe the evolution of the VMS compliance rate by modeling drivers' learning processes. We endogenize traffic dynamics, route and departure time choices, travel delays, and VMS compliance, and thereby captur their interactions and interdependencies in a holistic manner. A case study in the west end of Glasgow is carried out to understand the impact of VMS has on road congestion and route choices in both the short and long run. Our main findings include an adverse effect of the VMS on the network performance in the long run (the ``rebound" effect), and existence of an equilibrium state where both traffic and VMS compliance are stabilized.
\end{abstract}

\section{Introduction}

We consider a dual-time-scale dynamic traffic assignment (DTA) model on a road network in the presence of {\it variable message signs} (VMS). In the proposed model a discrete-time, day-to-day model of traffic evolution is coupled to a continuous-time, within-day model with both route and departure time choices. A central focus of this paper is to model and understand how VMS, together with user compliance, influence en route choices and travel delays in the within-day time scale, and how drivers' travel experiences (including experiences with the VMS) in turn affect their route and departure time choices as well as compliance with the VMS on a day-to-day time scale. Our goal is to not only assess and understand the effectiveness of VMS as an {\it intelligent transport system} (ITS) intervention to reduce congestion in the short run, but also to investigate its potential ``rebound" effect that could  lead to an adverse outcome in the long run. We propose to achieve this by endogenizing traffic dynamics, route and departure time choices, travel delays, and VMS compliance within our dual-time-scale DTA model and thereby analyzing their interactions and interdependencies in a holistic manner.

\subsection{Literature Review on VMS}

The impact of variable message signs (VMS) on traffic performance has been assessed in many studies
\citep{LC, MHMPJ, CHFB, CGMWZ, Wei}. Although the impact varies in different contexts (urban/extra urban, incident management/congestion mitigation), most of these studies conclude that variable message signs have a great potential to influence route choice, and hence traffic performance. However, an accurate estimate of such impact is strongly dependent on the level of compliance rate. A number of theoretical studies have assumed either full or partial compliance; in the latter case it is fixed or systematically varied to assess its effect on network performance \citep{CAT}. The key role played by the compliance rate on the assessment of the VMS impact has triggered significant attention from  researchers. One stream of research has focused on a realistic estimate of the compliance rate using different approaches: observations from traffic counts and loop detectors, surveys and virtual driving simulators \citep{RL, CHFB, Hoye, KHTS}. Another stream of research has investigated the key factors influencing the compliance rate \citep{Bonsall, VAKJY, WBS, Adler, CJ, BS}. These factors include: driver characteristics such as age, gender, familiarity with the network, information accuracy, clarity of the information provided (e.g. specification of the cause of the delay), visible queue and VMS location.

A few studies model the compliance rate as an endogenous variable, which may depend on a number of factors such as the ones mentioned above \citep{BS, CAT, DM, VAKJY}. Most of these studies recognize the key roles played by the information provided and previous experience in drivers' learning processes, and thereby developing models for route choice behavior that account for such learning processes and estimate the compliance rate accordingly. \cite{DM} adopts a fuzzy approach to estimate the perceived travel time resulting from the data fusion of two fuzzy variables: experience and  information provided. Their model accounts for the uncertainty associated with both variables. However, the model assumes no direct interaction among each driver's route choices.  \cite{CAT} propose a scenario-based probabilistic model to capture drivers' expected utilities associated with each of the alternatives (routes), and subsequently, their compliance with the advice. The impact of the compliance rate on network traffic and journey times, however, is not considered. A good attempt to capture the compliance rate and its dynamic variability is made by \cite{VAKJY} through a human-computer interactive experiment. The authors develop a sequential binary Logit choice model where the utilities are updated to reflect drivers' learning process. The two key variables included in the utility function are the perceived delay and perceived accuracy of information. An updating function is introduced for both variables to capture the learning process of drivers. Again, traffic delays are assigned randomly and independent of drivers' route choices.  \cite{BS} capture the variability of the compliance rate based on learning from experience. However, their study does not consider the effect route choices have on congestion and each individual driver's delay.

None of the abovementioned studies treat drivers' route choices influenced by travel information and travel delay simultaneously as endogenous quantities; as a result, the interaction between network congestion and driving behavior is insufficiently captured in these models. An exception is \cite{YY}, who  include the route choice behavioral model and consequent compliance rate estimation within a stochastic user equilibrium approach. There the compliance rate of drivers equipped with route guidance is calculated as the probability of the actual travel times being less than those of the non-compliant drivers. This model solves an equilibrium state and the compliance rate associated therein simultaneously without considering the evolutionary dynamics of the compliance rate or drivers' learning process.


\subsection{Contribution of This Paper}
Our literature review on VMS reveals a potential gap in simultaneously modeling the traffic dynamics and drivers' learning processes in an analytical (that is, non-simulation and non-experimental) framework, and in capturing their interaction and interdependencies in both short and long term. This paper aims to contribute to this line of research by proposing a dual-time-scale and analytical model that endogenizes  traffic dynamics, route and departure time choices, travel delays, and VMS compliance, and capture their interactions and interdependencies in a holistic manner. Our specific contributions are listed below.

\begin{itemize}
\item The within-day traffic dynamic is captured by a novel integration of the Lighthill-Whitham-Richards model \citep{LW, Richards} with en route updates influenced by the VMS. We propose an extension of the {\it dynamic network loading} (DNL) procedure, which is capable of calculating path travel times with a given set of path departure rates while modeling en route updates at any given location in the network with any user compliance rate. Moreover, such a DNL procedure captures realistic features of a dynamic traffic network such as shock waves and queue spillbacks.

\item We propose a day-to-day dynamic traffic assignment model with both route choices and departure time choices.  This model is new in that it incorporates drivers' bounded rationality when they adjust their travel arrangements.  

\item We propose two models for the day-to-day evolution of drivers' compliance rate with the VMS. These models assume that drivers' perception of the VMS is dependent on their past experience with the VMS and their familiarity with the network. Extensions of these two models that incorporates travel time variability and users' bounded rationality are also presented.

\item We apply the proposed methodology to a real-world test site located in Glasgow, Scotland. Various numerical studies are performed to investigate: (1) sensitivity of the model to its parameters; (2) short-term and long-term impact on route choices and user compliance; and (3) potential steady states of the traffic-compliance system and dynamic user equilibrium with en route path choices. 
\end{itemize}

\subsection{Organization} The rest of this paper is organized as follows. Section \ref{secNEB} serves as an introduction to basic mathematical notations and concepts employed by this paper. It also contains a brief description of the dynamic network loading (DNL) procedure. Section \ref{seccrdynamic} presents an analytical integration of the DNL procedure and the en route update influenced by the VMS. This is followed by several models that capture the evolution of user compliance with the VMS. Section \ref{secDTD} proposes a day-to-day DTA model that incorporates bounded user rationality. Finally, Section \ref{secCaseStudy} describes a case study in Glasgow.

\section{Notation and Essential Background}\label{secNEB}

Throughout this paper, we let $\tau\in\{1,\,2,\,3,\,\ldots,\}$ be one typical discrete day. The continuous time within each day is denoted by $t\in[t_0,\,t_f]$, where $[t_0,\,t_f]$ is a fixed commuting period, say, from 7 am to 10 am of each day.  We let $\mathcal{P}$ be the set of paths employed by travelers. For each path $p\in\mathcal{P}$ we define the path departure rate on day $\tau$, which is a function of departure time $t\in[t_0,\,t_f]$:
$$
h_p^{\tau}(\cdot):~[t_0,\,t_f]~\rightarrow~\Re_+
$$ 
where $\Re_+$ denotes the set of nonnegative real numbers. Each path departure rate $h^{\tau}_p(\cdot)$ is interpreted as the flow of departing vehicles measured at the entrance of the first arc of the relevant path.  We next define $h^{\tau}(\cdot)=\{h^{\tau}_p(\cdot): p\in\mathcal{P}\}$ to be a vector of departure rates. $h^{\tau}(\cdot)$ can be viewed as a vector-valued function of $t$, the departure time.

We denote the space of square integrable functions on the real interval $[t_0,\,t_f]$ by $L^2[t_0,\,t_f]$, and its subset consisting of non-negative functions by $L_+^2[t_0,\,t_f]$.  We stipulate that each path departure rate is square integrable; that is
$$
h^{\tau}_p(\cdot)\in L_+^2[t_0,\,t_f],\qquad \qquad h^{\tau}(\cdot)\in\big(L_+^2[t_0,\,t_f]\big)^{|\mathcal{P}|}
$$
where $\big(L_+^2[t_0,\,t_f]\big)^{|\mathcal{P}|}$ is the positive cone of  the $|\mathcal{P}|$-fold product of the Hilbert space $L^2[t_0,\,t_f]$.

Here, as in almost all dynamic traffic assignment modeling, the single most crucial ingredient is the network performance module encapsulated by the {\it path delay operator},  which maps a given vector of departure rates $h^{\tau}$ to a vector of path travel times. Mathematically, we let
\begin{equation*}
D_{p}(t,\,h^{\tau})\qquad \forall t\in[t_0,\,t_f],\quad  \forall p\in \mathcal{P}
\end{equation*}
be the path travel time of a driver departing at time $t$ and following path $p$, given the departure rate vector  $h^{\tau}$ associated with all the paths in the network on day $\tau$. 

We next define a generalized notion of travel cost by including not only path travel times, but also arrival penalties. The travel cost is defined as
\begin{equation}\label{cost}
\Psi _{p}(t,\,h^{\tau})~=~D_{p}(t,\,h^{\tau})+f\big( t+D_{p}(t,h^{\tau})-T_{A}\big) \qquad \forall t\in[t_0,\,t_f],\quad \forall p\in \mathcal{P}
\end{equation}
where $T_{A}$ is the desired arrival time which may depend on each O-D, and $T_{A}<t_{f}$. The term $f\big( t+D_{p}(t,h^{\tau})-T_{A}\big)$ assesses a nonnegative arrival penalty. We interpret $\Psi_p(t,\,h^{\tau})$ as the perceived travel cost of drivers departing at time $t$ following path $p$ given the vector of path departure rates $h^{\tau}$ on day $\tau$.

We write the demand satisfaction constraints as
\begin{equation}\label{cons}
\sum_{p\in \mathcal{P}_{ij}}\int_{t_0}^{t_f}h_p(t)\,dt~=~Q_{ij}\qquad\forall (i,\,j)\in\mathcal{W}
\end{equation}
where $\mathcal{P}_{ij}$ denotes the set of paths connection O-D  $(i,\,j)$, and $Q_{ij}$ is the (fixed) demand between $(i,\,j)$.  Using the notation
and concepts we have mentioned, the feasible region for the departure rate vector $h^{\tau}$ is
\begin{equation}\label{chapVI:lambda}
\Lambda~=~\left\{ h^{\tau}\geq 0:\sum_{p\in \mathcal{P}_{ij}}
\int_{t_{0}}^{t_{f}}h_{p}\left( t\right) dt=Q_{ij}\qquad \forall
\left( i,j\right) \in \mathcal{W}\right\} \subseteq \left( L_{+}^{2}[t_{0},\,t_{f}] \right) ^{\left\vert \mathcal{P}\right\vert }
\end{equation}

\subsection{The Dynamic Network Loading Procedure}\label{subsecDNL}
Our within-day dynamic is based on the DNL model formulated as a system of {\it differential algebraic equations} (DAEs) proposed by \citep{spillback}. This model employs the classical LWR link dynamics and captures vehicle spillback. Due to space limitation we omit details of this DNL model and refer the reader to \cite{spillback}.  In the next section we propose an analytical formulation of the en route update mechanism induced by the VMS by revising the vehicle turning ratios at a junction. Notice that such a formulation is fully generalizable to work with the cell transmission model \citep{CTM2} and the link transmission model \citep{LTM}.

\section{En Route Diversion and Day-to-Day Evolution of VMS Compliance}\label{seccrdynamic}

The first part of this section focuses on modeling en route update of drivers' path choices influenced by the VMS in a way consistent with the LWR model or any other physical queue models.  The second part of this section proposes several models that describes the dynamic evolution of drivers' compliance rates based on various assumptions.

\begin{figure}[h!]
\centering
\includegraphics[width=.8\textwidth]{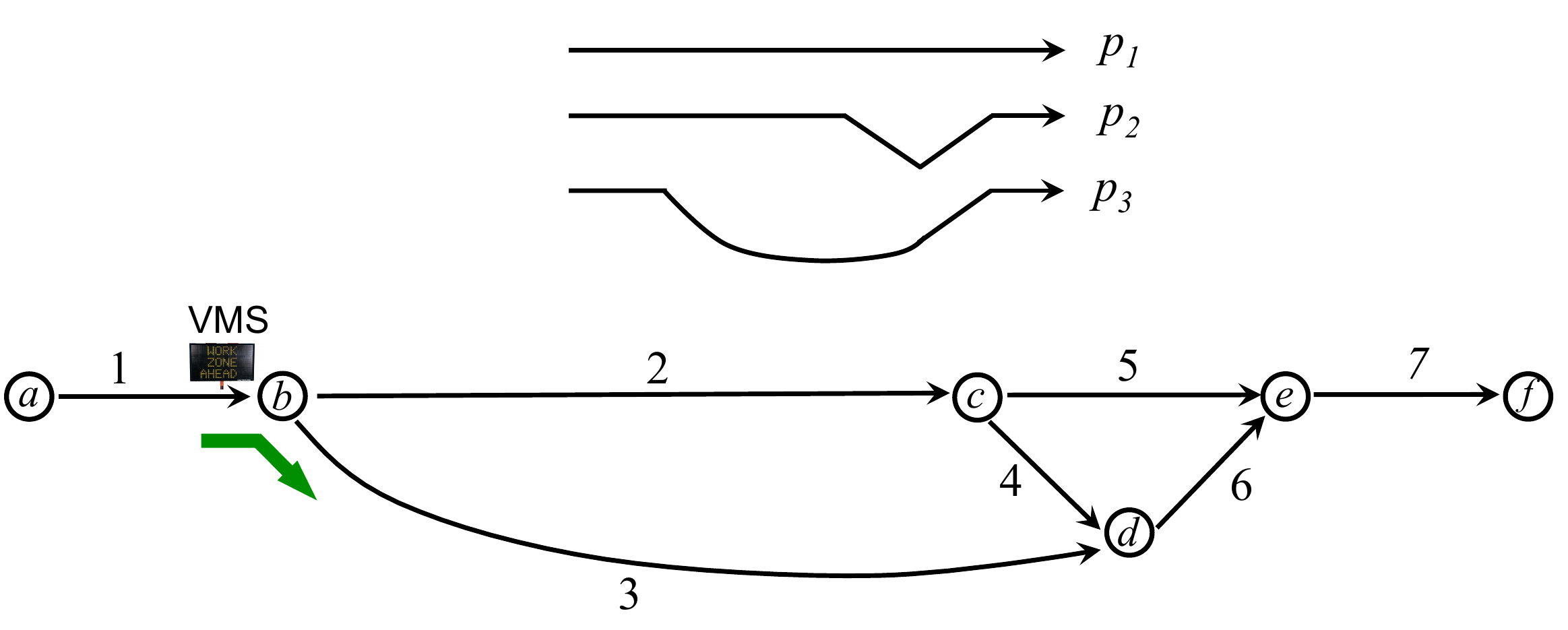}
\caption{A simple network with VMS.}
\label{figFig1}
\end{figure}

\subsection{Diversion of Flows at a Junction}
 Without loss of generality, we use a simple network shown in Figure \ref{figFig1} to illustrate the main idea for modeling en route diversion triggered by the VMS. The network contains three routes, $p_1$  $p_2$, $p_3$, and one origin-destination pair $(a,\,f)$. The VMS is located near the exit of link $1$. It is assumed that the VMS is always advising drivers to take route $p_3$, although such an assumption is easy to relax. Moreover, we allow the message sign to be on within time periods whose union is denoted by $\Omega\subset[t_0,\,t_f]$. For each day $\tau$, where $\tau=1,\,2,\,\ldots$, we let $h^{\tau}=\big(h^{\tau}_p(t),\,p\in\mathcal{P}\big)$ be the vector of path departure rate reviewed previously, where $\mathcal{P}=\{p_1,\,p_2,\,p_3\}$.

Throughout this paper, the VMS compliance rate is defined to be {\it the proportion of drivers who abandon their original routes in order to comply with what the VMS suggests}. If a driver is taking the VMS-suggested route as his/her route choice upon departure, then he/she is not influenced by the VMS and thus will not be accounted for by the compliance rate. 

 In order to capture route diversion induced by the VMS, we proceed as follows.  Denote by $\alpha_{1,2}(t)$ and $\alpha_{1,3}(t)$ the (time-varying) turning ratios at junction $b$ for downstream links $2$ and $3$ respectively. These quantities are determined within the dynamic network loading procedure based on path information, flow propagation constraints, and the first-in-first-out principle \citep{spillback}.  By definition, for every unit of flow exiting link $1$, $\alpha_{1,2}(t)$ of it is intended to advance into link $2$, thus if the VMS compliance rate on day $\tau$ is denoted by $CR^{\tau}$, the revised turning ratio, $\tilde\alpha_{1,2}(t)$, becomes 
\begin{equation}\label{revisedalpha1}
\tilde\alpha_{1,2}(t)~=~\alpha_{1,2}(t)-CR^{\tau}\alpha_{1,2}(t)
\end{equation}
Similarly, the revised turning ratio, $\tilde\alpha_{1,3}$, becomes
\begin{equation}\label{revisedalpha2}
\tilde\alpha_{1,3}(t)~=~\alpha_{1,3}(t)+CR^{\tau}\alpha_{1,2}(t)
\end{equation}
The revised turning ratios take into account route diversion suggested by the VMS with a fixed compliance rate. Using identity $\alpha_{1,2}(t)+\alpha_{1,3}(t)\equiv 1$ and the fact that $CR^{\tau}\in[0,\,1]$, one can easily verify that both revised turning ratios are between zero and one, and their sum equals one. Therefore, these revised turning ratios can be readily integrated into any junction model, which is embedded in a complete DNL procedure.

\subsection{Modeling Drivers'  Perception of the En Route Guidance and Their Compliance Rate}

Once the revised turning ratios are calculated according to \eqref{revisedalpha1} and \eqref{revisedalpha2}, we proceed as usual to complete the dynamic network loading (DNL) procedure with $\alpha_{1,2}(t)$ and $\alpha_{1,3}(t)$ replaced by $\tilde\alpha_{1,2}(t)$ and $\tilde\alpha_{1,3}(t)$ respectively. The DNL provides, as its output, the path travel times $D_p(t,\,h^{\tau})$, $p=p_1,\,p_2,\,p_3$, where $t$ denotes departure time from the origin of the path and $h^{\tau}$ denotes the vector of path departure rates on day $\tau$. The DNL procedure reported in  \cite{spillback} also calculates link exit times $\mu^{\tau}_a(t)$ for each link $a\in\mathcal{A}$ and  link entry time $t$. As usual, the superscript $\tau$ indicates a particular day.

It is straightforward to express the travel time of drivers who exited link $1$ at time $t$ and did not follow the guidance of the VMS as
$$
\mu_2^{\tau}\circ \mu_5^{\tau} \circ \mu_7^{\tau}(t) ~~\hbox{for path } p_1, \qquad \hbox{and}\qquad  \mu_2^{\tau}\circ\mu_4^{\tau}\circ \mu_6^{\tau} \circ \mu_7^{\tau}(t) ~~\hbox{for path } p_2
$$
where $\circ$ denotes composition of two functions; that is, $f\circ g(x)\doteq g\big(f(x)\big)$. The travel time of drivers who exited link $1$ at time $t$ and followed the guidance of the VMS is 
$$
\mu_3^{\tau}\circ\mu_6^{\tau}\circ \mu_7^{\tau}(t)
$$
Without causing any confusion, we employ shorter notations, $\mu^{\tau}_{2,5,7}(t)$, $\mu^{\tau}_{2,4,6,7}(t)$ and $\mu^{\tau}_{3,6,7}(t)$, for the three above-defined functions.

The compliance rate (CR) of a VMS is influenced by a number of factors such as those mentioned in the literature review, among which we will mainly focus on drivers' past experiences with such a route guidance, and drivers' own assessment of the suggested routes. Two additional extensions are also discussed to take into account drivers' bounded rationality and travel time variability. We make note of the fact that the proposed models are rather simplified; their purpose is solely to illustrate the behavior of the dual-time-scale traffic system with integrated traffic dynamics, route and departure time choices, and users' learning processes. We wish to analyze these models to provide more insights on the impact of the VMS on congestion and user behavior.

\subsubsection{Model I}

In this model, we assume that  drivers' perception of the effectiveness of the route guidance is directly related to the differences in the experienced (actual) journey times of drivers who followed such guidance and those who did not. In our particular example, such a difference may be expressed as 
\begin{equation}\label{examplesaving}
\mathcal{S}^{\tau}_{2\to 3}(t)~\doteq~{1\over 2} \big[\mu_{2,5,7}^{\tau}(t)+\mu_{2,4,6,7}^{\tau}(t)\big] -\mu_{3,6,7}^{\tau}(t)
\end{equation}
where the subscript $2\to 3$ indicates choosing link $3$ over link $2$. Notice that we simply average the travel times along the remaining portions of paths $p_1$ and $p_2$ to account for the journey times for non-compliant drivers, while noting that other more sophisticated treatments can also be applied without affecting our formulation. \eqref{examplesaving} expresses the direct saving in travel times by following the route guidance, which can be positive or negative.

We further aggregate the saving by averaging $\mathcal{S}^{\tau}_{2\to 3}(t)$ over $t$; the averaged saving is 
\begin{equation}\label{model1eqn1}
\bar{\mathcal{S}}_{2\to 3}^{\tau}~\doteq~{1\over |\Omega|}\int_{\Omega}  \mathcal{S}_{2\to3}^{\tau}(t)\,dt
\end{equation}
where we define $\Omega$ to be the union of time intervals during which the VMS is on. Eqn \eqref{model1eqn1} eliminates the time-variability of the saving and proposes a single (and possibly more robust) value as an indicator of the effectiveness of the route guidance. 

We then define drivers' perception of the route guidance as the expected  saving  in journey time by following the route guidance on day $\tau$. Mathematically, the perception on day $\tau$, denoted by $X^{\tau}$, is given by a weighted sum of the perception on a previous day, and the actual saving on day $\tau$:
\begin{equation}\label{model1eqn2}
X^{\tau}~=~(1-w) X^{\tau-1}+w \bar{\mathcal{S}}_{2\to 3}^{\tau} \qquad \tau~=~1,\,2\,\,\ldots
\end{equation}
for some weighting parameter $w\in(0,\,1)$. The updating rule of the form \eqref{model1eqn2} has been considered elsewhere too, for example, by \cite{VAKJY}. 

In deriving the compliance rate (CR), we consider two discrete choices: follow/not follow the route guidance displayed on the VMS. According to \eqref{examplesaving} and \eqref{model1eqn2}, the  drivers' perception (expected saving) of not following the guidance is clearly given as $-X^{\tau}$ on day $\tau$. Then by applying a logit model we get 
\begin{equation}\label{model1eqn3}
CR^{\tau+1}~=~{\exp(\beta X^{\tau})\over \exp(-\beta X^{\tau})+\exp(\beta X^{\tau})}
\end{equation}
where $\beta>0$.

\subsubsection{Model II}
A critical assumption made in the previous model is that drivers have perfect knowledge of the travel times on paths that he/she is currently not taking. This is reflected in the expressions for the savings $\mathcal{S}^{\tau}_{2\to3}(t)$ or $\bar{\mathcal{S}}_{2\to3}^{\tau}$. In the second model this assumption will be relaxed by assuming that drivers, when faced with an en route guidance, will make routing decisions based on their past experience with various paths in the network. Another important reason for considering such a model is the observation that the compliance rate with the VMS is affected  not only by their past experience with the VMS, but also by how the guidance provided in real-time is corroborated by the drivers' daily driving experiences \citep{Bonsall}.

To account for drivers' past experiences with various routes in the network, we let $\mu_{F}^{\tau}(t)$ to be the travel time starting from node $b$ (since this is where potential route diversion occurs) along path $p_3$, where ``$F$" stands for ``follow" (the VMS). In addition, let $\mu_{NF}^{\tau}(t)$ to be the path traversal time starting from $b$ along paths $p_1$ or $p_2$, where ``$NF$" stands for ``not follow" (the VMS). As usual, the superscript $\tau$ indicate a particular day and $t$ denotes the time at which the driver leaves node $b$. Mathematically, we have that 
\begin{equation}
\mu_{F}^{\tau}(t)~=~\mu^{\tau}_{3,6,7}(t),\qquad \mu_{NF}^{\tau}(t)~=~{1\over 2}\left(\mu_{2,5,7}^{\tau}(t)+\mu_{2,4,6,7}^{\tau}(t)\right)
\end{equation}
We eliminate the variability of path travel times in the within-day scale by averaging over $t$ to get:
\begin{equation}\label{avegeout}
\bar{\mu}_F^{\tau}~\doteq~{1\over t_f-t_0}\int_{t_0}^{t_f}\mu_{F}^{\tau}(t)\,dt,\qquad\qquad \bar \mu_{NF}^{\tau}~\doteq~{1\over t_f-t_0}\int_{t_0}^{t_f}\mu_{NF}^{\tau}(t)\,dt
\end{equation}
where $[t_0,\,t_f]$ represents the entire within-day assignment horizon under consideration. Similar to \eqref{model1eqn2}, we define the perceived path travel times $Y_{F}^{\tau}$ and $Y_{NF}^{\tau}$ respectively as
\begin{equation}\label{YFupdate}
Y_{F}^{\tau}~=~(1-w)Y_{F}^{\tau-1}+w \bar{\mu}_F^{\tau},\qquad Y_{NF}^{\tau}~=~(1-w)Y_{NF}^{\tau-1}+w \bar{\mu}_{NF}^{\tau},\qquad \tau~=~1,\,2\,\ldots
\end{equation}
for some weighting parameters $w\in(0,\,1)$. 

When experienced drivers face the en route guidance on day $\tau+1$, their macroscopic compliance rate, captured by the logit model based on their expected disutilities $Y_{F}^{\tau}$ and $Y_{NF}^{\tau}$, is give as
\begin{equation}
CR^{\tau+1}~=~{\exp(-\beta Y_{F}^{\tau})\over \exp(-\beta Y_{F}^{\tau})+\exp(-\beta Y_{NF}^{\tau})}
\end{equation}
where $\beta>0$.

\subsection{Model Extensions}

Models I and II proposed previously can be extended in a number of ways to capture more realistic decision-making processes. Due to space limitation, we proposes only two extension here, one for each model. The first extension incorporates a bounded user rationality; while the second one takes into account variability in travel times.


\subsubsection{Extension of Model I}
We postulate a threshold $\gamma>0$ such that drivers will only consider the VMS effective if it leads to a saving in travel time above $\gamma$. Thus, the average saving $\bar{\mathcal{S}}_{2\to 3}^{\tau}$ defined in \eqref{model1eqn1} can be revised to be 
\begin{equation}\label{revisedbarS}
\bar{\mathcal{S}}_{2\to3, \gamma}^{\tau}~=~\begin{cases}
0  \qquad   & \hbox{if} ~ \bar{\mathcal{S}}_{2\to3}^{\tau}\in[0,\,\gamma)
\\
\bar{\mathcal{S}}_{2\to3}^{\tau} \qquad  & \hbox{otherwise}
\end{cases}
\end{equation}
Notice that we use the subscript $\gamma$ to indicate the dependence on it. The rest of the model is the same as before and we call the extended model {\bf Model III}. The parameter $\gamma$ characterizes drivers' bounded rationality as it postulates certain indifference of drivers towards a positive yet insignificant gain by following the VMS, possibly due to route preference and imperfect information.

\subsubsection{Extension of Model II}
Both Model I and Model II proposed above disregard the travel time variability by averaging the travel time functions over time; see \eqref{model1eqn1} or \eqref{avegeout}. In reality the variation of travel time has a great impact on drivers' route choices, and it may outweigh the mean travel time in some circumstances \citep{Bonsall}.  Such a consideration can be easily incorporated into Model II by defining 
\begin{equation}\label{moment}
\sigma_{F}^{\tau}~\doteq~\left( {1\over t_f-t_0}\int_{t_0}^{t_f}(\mu_{F}^{\tau}(t)-\bar\mu_F^{\tau})^2\,dt\right)^{{1\over 2}},\qquad \sigma_{NF}^{\tau}~\doteq~\left( {1\over t_f-t_0}\int_{t_0}^{t_f}(\mu_{NF}^{\tau}(t)-\bar\mu_{NF}^{\tau})^2\,dt\right)^{{1\over 2}}
\end{equation}
and writing 
\begin{equation}\label{YFupdatenew}
Y_{F}^{\tau}~=~(1-w)Y_{F}^{\tau-1}+w\cdot \bar{\mu}_F^{\tau}\cdot \sigma_{F}^{\tau},\qquad Y_{NF}^{\tau}~=~(1-w)Y_{NF}^{\tau-1}+w\cdot  \bar{\mu}_{NF}^{\tau}\cdot \sigma_{NF}^{\tau},\qquad \tau~=~1,\,2\,\ldots
\end{equation}
Notice that in discrete time,  $\bar\mu_F^{\tau}$ and $\sigma_F^{\tau}$ represents mean and standard deviation of the travel time vector, respectively. The rest of the model remains the same. We will call such an extension {\bf Model IV}. 

\section{A Day-to-Day Assignment Model Incorporating Bounded User Rationality}\label{secDTD}
This section discusses the day-to-day component of the model that describes the dynamic evolution of drivers' travel decisions. It is assumed that drivers adjust their departure time and route choices on a daily basis in search for a more efficient travel arrangement. 

We propose a day-to-day dynamic traffic assignment model based on {\it boundedly rational dynamic user equilibrium} (BR-DUE) proposed by \cite{BRDUE} and the fixed-point iteration algorithm associated therein. The notion of bounded rationality is a relaxation of the Wardropian equilibrium assumption which is based on perfect user rationality.  That is, the BR postulates a cost tolerance,  also known as {\it indifference band} \citep{MC},  such that drivers will only abandon their current travel choices if an alternative offers cost saving that is beyond such a tolerance. 

We let $\varepsilon=(\varepsilon_{ij}^{p},\,p\in\mathcal{P}_{ij},\,(i,\,j)\in\mathcal{W})$ be the vector of (non-negative) tolerances, each of which depends on the O-D $(i,\,j)$ and the path $p$, where $\mathcal{P}_{ij}$ denotes the set of paths connecting $(i,\,j)$, and $\mathcal{W}$ is the set of O-D pairs. Notice that unlike many existing studies on BR which assume that the tolerance depends only on the O-D, we allow it to depend also on the path. This additional relaxation can be applied to model drivers' preferences towards a particular set of paths. 

The fixed-point update proposed by \cite{BRDUE} can be mathematically expressed as
\begin{equation}\label{fpu}
h^{\tau+1}(\cdot)~=~P_{\Lambda}\big[h^{\tau}(\cdot) - \lambda\Phi^{\varepsilon}(\cdot,\,h^{\tau})\big]\qquad\tau~=~1,\,2,\,\ldots
\end{equation}
where $\lambda>0$ is a fixed constant. The operator $P_{\Lambda}$ means the minimum-norm projection onto the convex set $\Lambda$ defined in \eqref{chapVI:lambda}. Moreover, $\Phi^{\varepsilon}(t,\,h^{\tau})=\big(\Phi_p^{\varepsilon}(t,\,h^{\tau}),\, p\in\mathcal{P}\big)$ and 
\begin{equation}\label{Phidef1}
\Phi^{\varepsilon}_p(t,\,h^{\tau})~=~\max\left\{\Psi_p(t,\,h^{\tau}),~ v_{ij}(h^{\tau})+\varepsilon_{ij}^p\right\}-\left(\varepsilon_{ij}^p-\min_{q\in\mathcal{P}_{ij}}\left\{\varepsilon_{ij}^q\right\}\right)\qquad \forall p\in\mathcal{P}_{ij}
\end{equation}
where $v_{ij}(h^{\tau})$ denotes the minimum travel cost within $(i,\,j)$. \eqref{fpu} can be explicitly instantiated, leading to the following day-to-day DTA model (we omit details of the derivation and refer the reader to \cite{BRDUE}): 

\vspace*{10pt}

\noindent \noindent \fbox{Step 0. Initialization.} Set day $\tau=1$ and identify an initial feasible $h^{\tau}\in\Lambda$.

\vspace*{10pt}

\noindent \noindent \fbox{Step 1. Dynamic network loading.}  Solve the dynamic network loading problem with $h^{\tau}$ and compliance rate $CR^{\tau}$ to obtain $\Psi_p(t,\,h^k)$, $p\in\mathcal{P},\,t\in[t_0,\,t_f]$. Then find the following quantities
$$
\Phi^{\varepsilon}_p(t,\,h^{\tau})~=~\max\left\{\Psi_p(t,\,h^{\tau}),~ v_{ij}(h^{\tau})+\varepsilon_{ij}^p\right\}-\big(\varepsilon_{ij}^p-\min_{q\in\mathcal{P}_{ij}}\big\{\varepsilon_{ij}^q\big\}\big)\qquad \forall p\in\mathcal{P}_{ij},\quad \forall t\in[t_0,\,t_f]
$$

\noindent \noindent \fbox{Step 2. Find the dual variable} For each $(i,\,j)\in\mathcal{W}$, solve the following nonlinear equation for $\eta_{ij}$, using root-search algorithms. 
$$
\sum_{p\in \mathcal{P}_{ij}}\int_{t_0}^{t_f}\left[ h_{p}^{k}(t)
-\lambda \Phi^{\varepsilon} _{p}\left( t,h^{k}\right) +\eta_{ij}\right] _{+}dt~=~Q_{ij}
$$
where $[u]_+\doteq \max\{0,\,u\}$.

\noindent \noindent \fbox{Step 3. Update path departure rates.} For each $(i,\,j)\in\mathcal{W}$ and $p\in\mathcal{P}_{ij}$, update the departure rate as
$$
h_{p}^{\tau+1}(t)~=~\left[ h_{p}^{\tau}(t)-\lambda \Phi^{\varepsilon} _{p}\left( t,h^{\tau}\right) +\eta_{ij}
\right] _{+}\qquad\forall t\in[t_0,\,t_f]
$$ 
Set $\tau=\tau+1$, repeat Step 1 through Step 3.

\section{Case Study}\label{secCaseStudy}

\subsection{Description of The Glasgow Test Site}

\begin{figure}[h!]
\centering
\includegraphics[width=.8\textwidth]{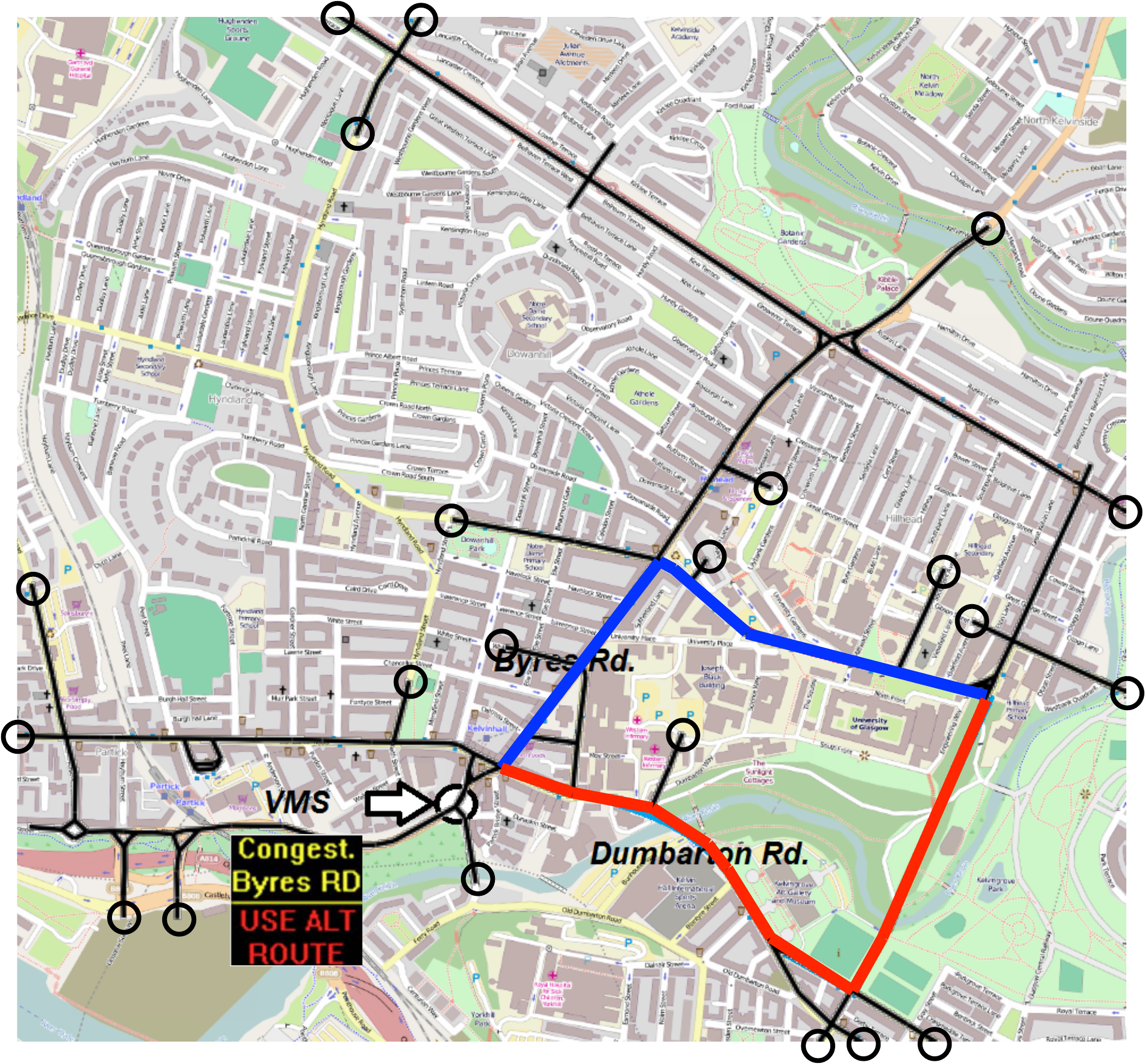}
\caption{The Glasgow test network highlighted in black. The circles represent the 21 traffic zones giving rise to 420 O-D pairs.  The location of the VMS is shown with arrows. Two main corridors are highlighted with blue and red; the latter is recommended by the VMS, as shown by the example message.}
\label{figVMSGlasgow}
\end{figure}

The west end of Glasgow has been chosen as our test network, see Figure \ref{figVMSGlasgow}. The test site has been identified by the local traffic authority and environmental agencies to suffer from traffic and air quality problems. Byres Road (highlighted in blue)  is often affected by severe congestion as it not only connects the radial routes to the city centre for drivers approaching Glasgow from the west, but also provides access to the university and other local destinations. The Dumbarton Road corridor (highlighted in red) serves as an alternative route for drivers traveling from the west end towards east (city centre). Our numerical study is based on preliminary work aiming at calibrating and validating traffic flows and counts obtained from loop data and surveys, and represents  morning peak hour (7:30-9:30) on Monday -Thursday. 

This case study is part of the CARBOTRAF project funded by the 7th Framework (http://www.carbotraf.eu). The use of VMS as an ITS intervention, among others such as coordinated signal control, serve the general purpose of reducing congestion and emissions of CO$_2$ and Black Carbon.

 There are 21 zones in the network leading to 420 OD-pairs. The demand matrix for the test network was derived from the Scotland O-D matrix and adjusted this OD matrix to the Glasgow test site by using traffic count data from several key network locations.  The network has 478 links; and a total of 481 paths are generated for the route-based dynamic traffic assignment model. Within this network 22 O-D pairs and 86 paths are affected by the VMS. Naturally, the VMS compliance rate and the learning process discussed in Section \ref{seccrdynamic} are carried out within each O-D pair. In other words, they are O-D specific.

\subsection{Numerical Results}

The numerical results are reported in three areas identified in the introduction: (1) model sensitivity to parameters; (2) impact on route choices; and (3) steady state and convergence to an equilibrium.

\subsubsection{Model Sensitivity to Parameters}
The proposed models employ a number of parameters that need to be calibrated using empirical data or surveys. As a preliminary step, we conduct an informal sensitivity analysis in this paper that highlights the qualitative behavior of the model when these parameters vary. Table \ref{tabsen} summarizes the model sensitivity to perturbations in these parameters (due to space limitation details of the sensitivity analyses are omitted). According to our findings, the models are most sensitive to $\beta$, which is the logit coefficient that reflects drivers' sensitivity to the costs (or difference in the costs). Moreover, $w$, which is the weighting parameter used to construct drivers' perception \eqref{model1eqn2} or \eqref{YFupdate}, may also influence the result to some extent. The parameter $\gamma$, which measures users' bounded rationality in Model I, has a similar impact on the models as $w$. All the models are  insensitive to the initial values such as $X^0$, $Y_F^0$ and $Y_{NF}^0$.

\begin{table}[h!]
\centering
\begin{tabular}{|c|c|c|c|c|}
\hline 
     & Model I    & Model II   & Model III   & Model IV
 \\\hline
$X^0$; $[ -600 ,\, 600]$  &   Low   &   -- &  Low    &  -- 
\\\hline  
$Y_F^0$ or $Y_{NF}^0$; $[100,\,600]$   &  --  &  Low  &  --  & Low
\\\hline
$w$; $[ 0.2 ,\, 0.5]$  & Medium        &  Low     &    Medium & Low
\\\hline
$\beta$; $[ 0.001 ,\, 0.1]$     &  High   &  High   & High   &  High
\\\hline
$\gamma$; $[0,\,300]$   & --  &  --   &  Medium   &  --
\\
\hline
\end{tabular}
\caption{Model sensitivity to the parameters. The interval following each parameter represents the range of values tested.}
\label{tabsen}
\end{table}

\begin{figure}[h!]
\centering
\includegraphics[width=\textwidth]{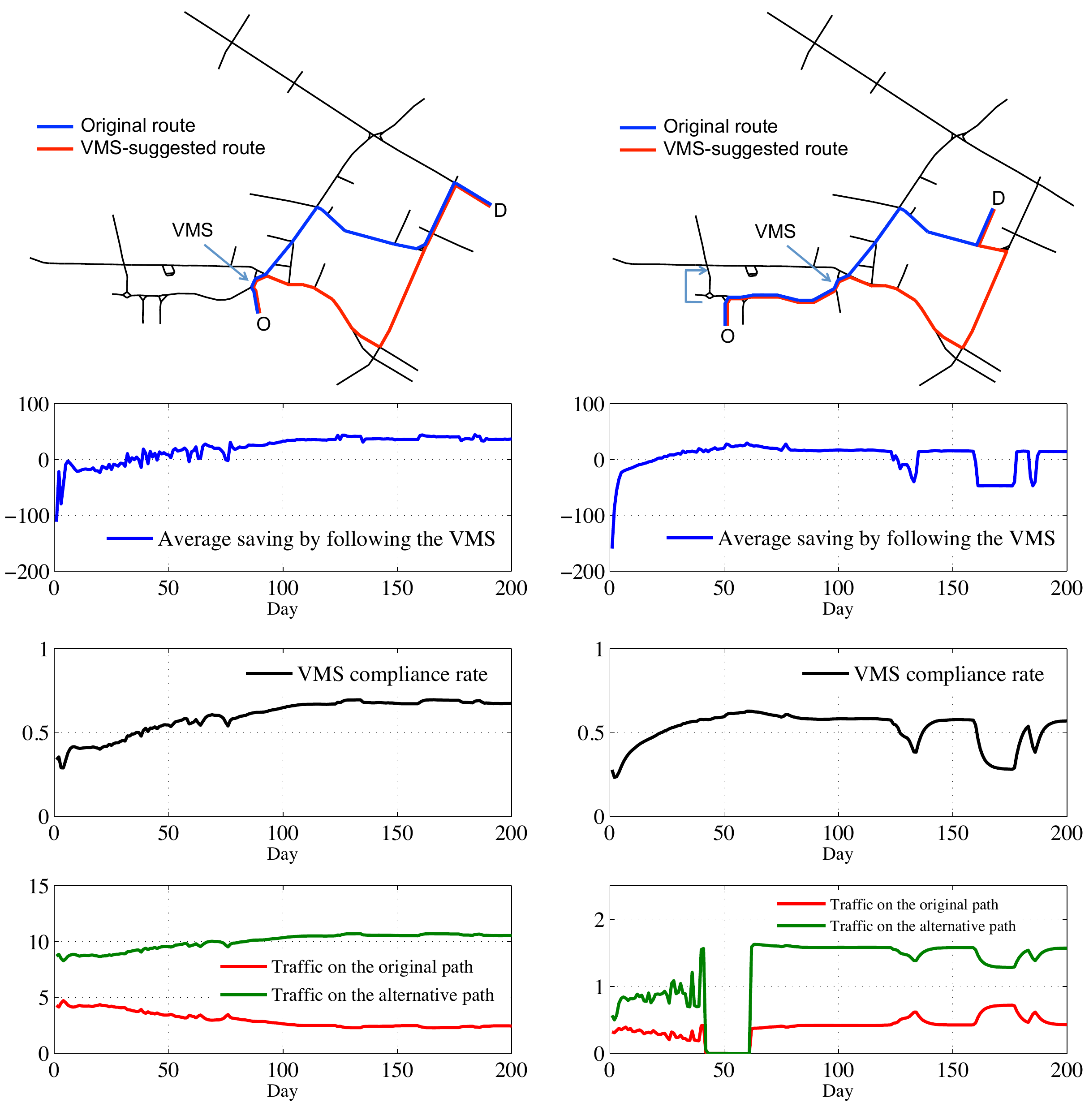}
\caption{Computational results of Model I for two O-D pairs (left and right column).}
\label{figTask7_term}
\end{figure}

\subsubsection{Impact on Route Choices}
We use a specific scenario to illustrate the short-term and long-term effect of VMS on route choices and user compliance. Consider Model I with parameters $w=0.3$, $\beta=0.01$. The computation is carried out for 200 days, or time periods. We illustrate the solution using two specific O-D pairs shown in the first row of Figure \ref{figTask7_term}. Two main routes are highlighted for each O-D pair; the original (primary) route is mark blue, and the alternative route suggested by the VMS is marked red. For the first O-D (corresponding to the first column of Figure \ref{figTask7_term}), the alternative route starts off as the less efficient route since it has a longer travel distance; this is seen by the negative saving during the first 40 days. Later on, however, as traffic evolves the alternative route becomes more efficient and actually saves travel time compared to the original route. This is caused by the increased congestion on the original corridor (highlighted blue). As a consequence, the compliance rate starts to climb (third row), and more traffic shift towards the alternative route (highlighted red); see the forth row of the figure. After about 100 days, the alternative route becomes saturated, causing the saving in travel time by following the VMS to stop increasing. Thus no more traffic shift from the original path to the alternative path; and the compliance rate also stabilizes. 

For the second O-D (right column of Figure \ref{figTask7_term}), the situation is similar to the first O-D during the first 100 days. (However, one complication exists for this O-D due to the presence of another set of paths that do not go through the VMS, which is indicated by an arrow in Figure \ref{figTask7_term}. This explains the period around the 50th day when no traffic exists on the two paths affected by the VMS; see the fourth row.) Unlike the previous case, the saving within the second O-D by following the VMS undergo three major drops to below zero during the period $[100, \, 200]$ (days), which leads the compliance rate to drop as well. This is caused by the increasing congestion on the alternative route, which is in turn due to drivers following the VMS (the rebound effect).

\subsubsection{Stationary State and Convergence}
\begin{figure}[h!]
\centering
\includegraphics[width=\textwidth]{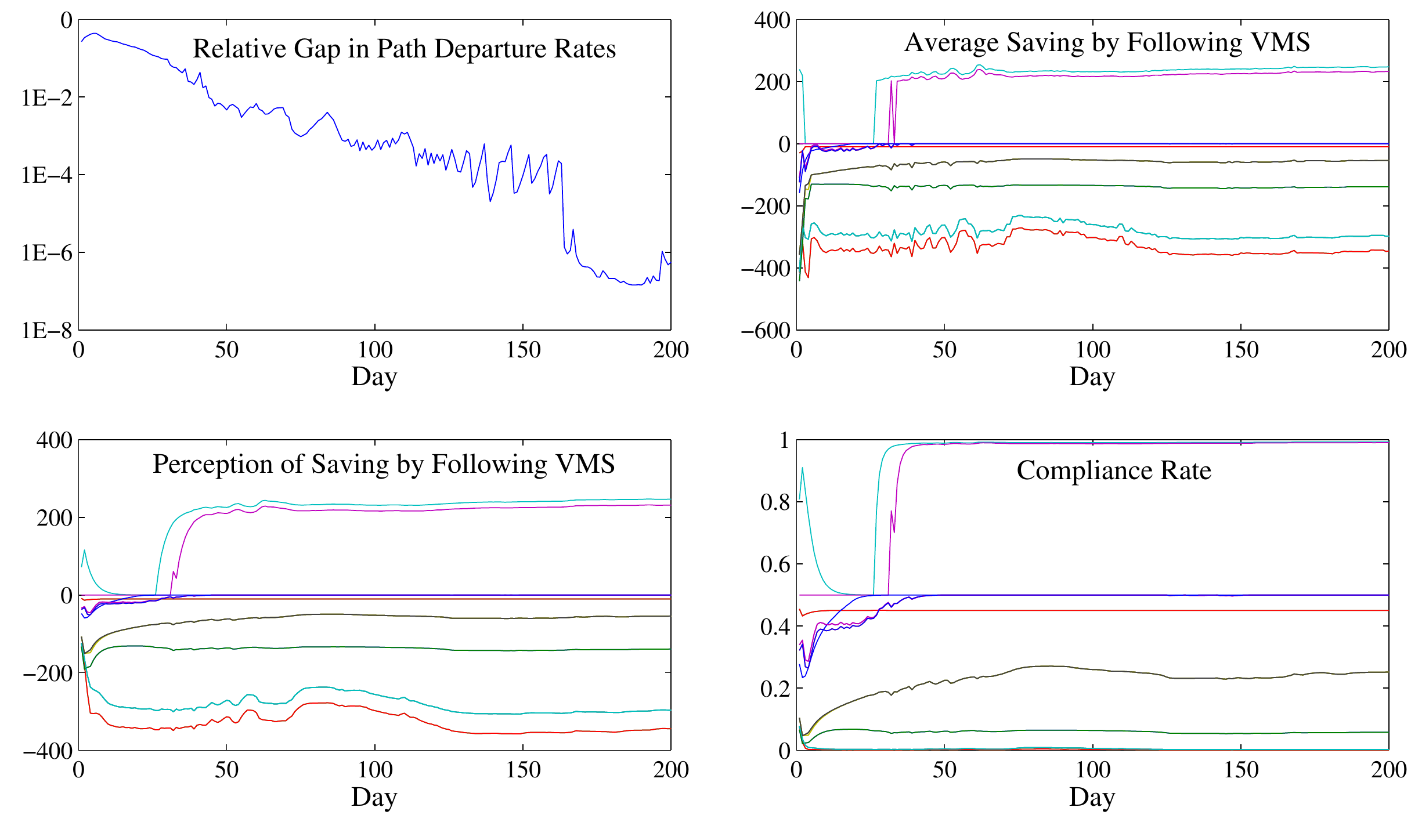}
\caption{Convergence of path departure rates (upper left), average saving by following VMS (upper right), drivers' perception of savings by following VMS (lower left), and the VMS compliance rate (lower right). The last three quantities are O-D specific; thus each color represents an O-D pair (there are totally 22 relevant O-D pairs).}
\label{figTask6_conv}
\end{figure}

Our proposed subject of study is a special case of the day-to-day dynamic traffic assignment models. Naturally, we would like to search for steady states of such a dynamical system. Figure \ref{figTask6_conv} provides confirmation of convergence to a steady state, when Model III with parameters $w=0.3,\,\beta=0.01,\,\gamma=200s$ is employed. From this figure we observe that both the path departure rates and the compliance rates stabilizes, as the number of ``days" reach 200. Notice that the convergence of the departure rates is measured by a relative gap defined by 
$$
\hbox{Relative gap in path departure rates}~=~{\|h^{\tau}-h^{\tau-1}\|_{L^2}\over \|h^{\tau-1}\|_{L^2}}
$$
where $\|\cdot\|_{L^2}$ is the norm on the Hilbert space $\big(L^2[t_0,\,t_f]\big)^{|\mathcal{P}|}$. In addition, the savings by following the VMS on a particular day, as well as the overall perception on such a saving, converge to constant values (see Figure \ref{figTask6_conv}). 

This example shows the existence of a boundedly rational dynamic user equilibrium with route and departure time choices, and in the presence of active route guidance in the network -- a type of user equilibrium not studied before. However, we did not observe convergence for all the models/parameters chosen (for example, the right column of Figure \ref{figTask7_term}). Further research is undertaken to analytically formulate such a class of user equilibria and investigate their existence.

\section{Conclusion and future research}
We propose a day-to-day dynamic traffic assignment model with en route travel updates influenced by the VMS. Traffic dynamics and user's learning processes are simultaneously modeled and their interactions and interdependencies analyzed. The results show very complex behavior of the traffic system when both route choices and VMS compliance are endogenized. Furthermore, we observe in many cases convergence of both traffic and compliance rate to a steady state, indicating the existence of a dynamic user equilibrium (DUE). Further work is under way to investigate theoretically existence result for such a DUE. Moreover, the complex behavior mentioned before is used as the basis to construct further insights and strategies for traffic management and control.

\bibliographystyle{model2-names}
\bibliography{<your-bib-database>}

\begin{thebibliography}{00}




\bibitem[Adler, 2001]{Adler} Adler, J. L., 2001. Investigating the learning effects of route guidance and traffic advisories on route choice behavior. Transportation Research Part C 9 (1), 1-14.



\bibitem[Ben-Elia and Shiftan, 2010]{BS} Ben-Elia, E., Shiftan, Y., 2010. Which road do I take? A learning-based model of route-choice behavior with real-time information. Transportation Research Part A 44 (4), 249-264.


	
\bibitem[Bonsall, 1992]{Bonsall} Bonsall, P., 1992. The influence of route guidance advice on route choice in urban networks. Transportation 19 (1), 1-23.
	
	
	

\bibitem[Chatterjee et al., 2002]{CHFB} Chatterjee, K., Hounsell, N. B., Firmin, P.E., Bonsall, P.W., 2002. Driver response to variable message sign information in London. Transportation Research Part C 10 (2), 149-169.


\bibitem[Chen et al., 2008]{CGMWZ} Chen S.,  Liu, M., Gao, L., Meng C., Li W., Zheng J., 2008. Effects of variable message signs (VMS) for improving congestions. International Workshop on Modelling, Simulation and Optimization, 416-419. 


\bibitem[Chen and Jovanis, 2003]{CJ} Chen, W., Jovanis, P., 2003. Analysis of driver en-route guidance compliance and driver with ATIS using a travel simulation experiment. 81st Transportation Research Board Annual Meeting, Washington, DC.


\bibitem[Chorus et al., 2009]{CAT} Chorus, C. G., Arentze, T. A.,  Timmermans, H.J.P., 2009. Traveler compliance with advice: A Bayesian utilitarian perspective. Transportation Research Part E 45 (3), 486-500.


\bibitem[Daganzo, 1995a]{CTM2}
Daganzo, C.F., 1995a. The cell transmission model. Part II: Network traffic. Transportation Research Part B 29 (2), 79-93.





\bibitem[Dell'Orco and Marinelli, 2009]{DM} DellÕOrco, M., Marinelli, M., 2009. Fuzzy data fusion for updating information in modeling driversÕ choice behavior. Emerging Intelligent Computing Technology and Applications. With Aspects of Artificial Intelligence. D.-S. Huang, K.-H. Jo, H.-H. Lee, H.-J. Kang and V. Bevilacqua, Springer Berlin Heidelberg. 5755: 1075-1084.

\bibitem[Friesz et al., 2013b]{FHNMY}
Friesz T.L., Han, K., Neto, P.A., Meimand, A., Yao, T., 2013b. Dynamic user equilibrium based on a hydrodynamic model. Transportation Research Part B 47 (1), 102-126.



\bibitem[Garavello and Piccoli, 2006]{GP}
Garavello, M.,  Piccoli,  B., 2006.  Traffic Flow on 
Networks. Conservation Laws Models.  AIMS Series on Applied Mathematics, Springfield, Mo..



\bibitem[Garavello and Piccoli, 2009]{Garavello and Piccoli 2009} Garavello, M., Piccoli, B., 2009. Conservation laws on complex networks. Ann. Inst. H. Poincar Anal. Non LinÃ©aire 26 (5),  1925-1951.




\bibitem[Han et al., 2014a]{spillback} Han, K., Friesz, T.L., Yao, T., 2014a. Vehicle spillback on dynamic traffic networks and what it means for dynamic traffic assignment models. 5th International Symposium on Dynamic Traffic Assignment. Salerno, Italy, 17-19 June 2014. 

\bibitem[Han et al., 2014b]{HGPFY} Han, K., Gayah, V., Piccoli, B., Friesz, T.L., Yao, T. 2014b. On the continuum approximation of the on-and-off signal control on dynamic traffic networks. Transportation Research Part B: Methodological, 61, 73-97.

\bibitem[Han et al., 2014c]{BRDUE} Han, K., Szeto, W.Y., Friesz, T.L., 2014c. Formulation, existence, and computation of simultaneous route-and-departure choice bounded rationality dynamic user equilibrium with fixed or endogenous user tolerance. Under review. Preprint available at: http://arxiv.org/abs/1402.1211


\bibitem[Hoye, 2011]{Hoye} Hoye, A., Sorensen, M., Elvik, R., Akhtar, J., Navestad, T., Vaa, T., 2011. Evaluation of variable message signs in Trondheim. N. C. f. T. Research, Oslo.




\bibitem[Kattan et al., 2011]{KHTS} Kattan, L., Habib, K.M.N., Tazul, I., Shahid, N., 2011. Information provision and driver compliance to advanced traveller information system application: case study on the interaction between variable message sign and other sources of traffic updates in Calgary, Canada. Canadian Journal of Civil Engineering 38 (12), 1335-1346.




\bibitem[Lam and Chan, 1996]{LC} Lam, W.H.K., Chan, K.S., 1996. A stochastic traffic assignment model for road network with travel time information via variable message signs. Intelligent Vehicles Symposium, 1996. Proceedings of the 1996 IEEE, 99-104.


\bibitem[Lebacque and Khoshyaran, 1999]{Lebacque and Khoshyaran 1999} Lebacque, J., Khoshyaran, M., 1999. Modeling vehicular traffic flow on networks using macroscopic models, in Finite Volumes for Complex Applications II, 551-558, Hermes Science Publications, Paris. 



\bibitem[Lighthill and Whitham, 1955]{LW} Lighthill, M.,   Whitham, G., 1955. On kinematic waves. II. A theory of traffic flow on long crowded roads.  Proceedings of the Royal Society of London. Series A, Mathematical and Physical Sciences 229 (1178),  317-345.

\bibitem[Mahmassani and Chang, 1987]{MC} Mahmassani, H., Chang, G., 1987. On boundedly rational user equilibrium in transportation systems. Transportation Science 21 (2), 89-99.


\bibitem[Mammar et al., 1996]{MHMPJ} Mammar, S., Hai-Salem, H., Messmer, A., Papageorgiou, M., Jensen, L., 1996. Vehicle Navigation and Information Systems Conference, 1996. VNIS '96, 12-22. 

\bibitem[Newell, 1993]{Newell}
Newell, G.F., 1993.  A simplified theory of kinematic waves in highway traffic, part I: General theory. Transportation Research Part B 27 (4), 281-287.



\bibitem[Ramsey and Luk, 1997]{RL} Ramsey, E.D., Luk, J., 1997. Route choice under two australian travel information systems, ARRB Transport Research.


\bibitem[Richards, 1956]{Richards} Richards, P.I., 1956. Shockwaves on the highway. Operations Research  4 (1), 42-51.



\bibitem[Vaughan et al., 1993]{VAKJY} Vaughn, K.M., Abdel-aty, M.A., Kitamura, R., Jovanis, P.P., Yang, H., 1993. Experimental Analysis And Modeling Of Sequential Route Choice Under ATIS In A Simple Traffic Network. Transportation Research Record 1408, 75-82.

\bibitem[Wardman et al., 1997]{WBS} Wardman, M., Bonsall, P.W., Shires, J.D., 1997. Driver response to variable message signs: a stated preference investigation. Transportation Research Part C 5 (6), 389-405.


\bibitem[Wei et al., 2009]{Wei} Wei, S., Wu, J., Zhou, S., Zhang, L., Tang, Z., Yin, Y.Y., Kuang, L., Wu, Z., 2009. Variable message sign and dynamic regional traffic guidance. Intelligent Transportation Systems Magazine, IEEE 1 (3), 15-21.

\bibitem[Wardrop, 1952]{Wardrop} Wardrop, J., 1952. Some theoretical aspects of road traffic research. In ICE Proceedings: Part II, Engineering Divisions 1, 325-362. 




\bibitem[Yin and Yang, 2003]{YY} Yin, Y., Yang, H., 2003. Simultaneous determination of the equilibrium market penetration and compliance rate of advanced traveler information systems. Transportation Research Part A  37 (2), 165-181.



\bibitem[Yperman et al., 2005]{LTM} Yperman, I., S. Logghe and L. Immers, 2005. The Link Transmission Model: An Efficient Implementation of the Kinematic Wave Theory in Traffic Networks, Advanced OR and AI Methods in Transportation, Proc. 10th EWGT Meeting and 16th Mini-EURO Conference, Poznan, Poland, 122-127, Publishing House of Poznan University of Technology.



 \end{thebibliography}



\end{document}